\documentclass[11pt]{amsart}
\usepackage{amsfonts,amsmath,amssymb}
\usepackage[all]{xy}
\usepackage{hyperref}
\usepackage[OT2,T1]{fontenc}

\newcommand\textcyr[1]{{\fontencoding{OT2}\fontfamily{wncyr}\selectfont #1}}

\begin{document}

\parindent=0pt
\parskip=6pt

\newcommand{\ie}{\textit{i}.\textit{e}.\,}
\newcommand{\eg}{\textit{e}.\textit{g}.\,}
\newcommand{\cf}{{\textit{cf}.\,}}

\newcommand{\CF}{{\sf CF}}
\newcommand{\aut}{{\rm Aut}}
\newcommand{\iN}{{\rm In}}
\newcommand{\B}{{\mathbb B}}
\newcommand{\C}{{\mathbb C}}
\newcommand{\R}{{\mathbb R}}
\newcommand{\mH}{{\mathbb H}}
\newcommand{\SO}{{\rm SO}}
\newcommand{\SU}{{\rm SU}}
\newcommand{\Z}{{\mathbb Z}}
\newcommand{\op}{{\rm op}}
\newcommand{\Sl}{{\rm Sl}}
\newcommand{\F}{{\mathbb F}}
\newcommand{\half}{{\textstyle{\frac{1}{2}}}}
\newcommand{\dslash}{{/\!\!/}}
\newcommand{\im}{{\rm im}}
\newcommand{\ve}{{\varepsilon}}

\title{On the canonical formula of C L\'evi-Strauss, II}

\author[Jack Morava]{Jack Morava}

\address{Department of Mathematics, The Johns Hopkins University,
Baltimore, Maryland 21218}

\email{jack@math.jhu.edu}

\date{Imbolc 2020}

\maketitle 

{\bf \S 1} {\sc Introduction and Organization}

\begin{quotation}{\noindent We venture a leap: we grant {\it ab initio} that
there is `something there' to be translated \dots \medskip

George Steiner, quoted in \cite{3} (p 19)}\end{quotation}

{\bf 1.1} This is a sequel to, and elaboration of, earlier work \cite{21,22} 
on a possible formal model for the canonical formula
\[
\CF  : F_x(a):F_y(b) \; \simeq \; F_x(b):F_{a^{-1}}(y)
\]
of C L\'evi-Strauss, which he proposed \cite{1,14,15,16}, \cf \: [Appendix] as 
a tool for the structural analysis of mythological systems; see \cite{12} (p 
562) for a discussion of the hair-raising example
\[
{\rm marriage}_{\rm solidarity} : {\rm rape}_{\rm hostility} \simeq
{\rm marriage}_{\rm hostility} : {\rm dissociation}_{\rm rape} \;.
\]
The model proposed here is based on the study of small finite groups, which
have proved useful in the classification of kinship systems, crystallography,
and other fields \cite{8,19,25,27}; a short account of some of the mathematics
involved is postponed till {\bf \S 3}. For clarity, a sketch of the
model, with technical details backgrounded, is displayed immediately below.
However, because translation across the conceptual and cognitive gulfs
separating anthropology and mathematics raises significant questions, 
a short discussion ({\bf \S 2}) of models and metaphors precedes the (not 
actually very complicated) verification of the details of the model, in 
{\bf \S 4}.

{\bf 1.2} {\sc A model} 

\begin{quotation}{\noindent The first formulation of the real situation
would be to state that the speaking subject  perceives neither idea $a$ nor 
form $A$, but only the relation $a/A$ \dots What he perceives is the relation
between the two relations $a/AHZ$ and $abc/A$, or $b/ARS$ and $blr/B$, etc. 
This is what we term the LAST QUATERNION \dots \medskip

F de Saussure, {\bf Writings}\dots \cite{26} (p 22), \:\cite {16}}
\end{quotation}

{\bf Proposition:} {\it To elements $\{x,a,y,b\}$, \eg $\{1,i,j,k\}$ of the
group $Q$ of Lipschitz unit quaternions [3.2.1], the function
\[
x,a \mapsto \Phi_x(a) = 2^{-1/2}(x - a) \in 2 \cdot O
\]
assigns elements of the binary octahedral group [3.4.2] such that the 
anti-automorphism $\lambda = *I\sigma^2$ [3.3{\it iii}] of that group 
transforms the (noncommutative [3.2.2]) ratio $\{\Phi_x(a):\Phi_y(b)\}$ into} 
\[
\{\Phi_{\lambda(x)}(\lambda(a)):\Phi_{\lambda(y)}(\lambda(b))\} \; = \;
\{\Phi_x(b):\Phi_{a^{-1}}(y)\} \;.
\]

Terms such as $x,a,y,b$ or $1,i,j,k$ will be referred to here as `values', 
while functions such as $F$ and $\Phi \in 2 \cdot O$ of ordered pairs of 
values will be called `valences'. [Citations such as [1.2] refer to sections 
of this paper.] \bigskip

{\bf \S 2} {\sc Wider Questions} \bigskip
\begin{quotation}{\noindent It takes a while to learn that things like
tenses and articles \dots are not `understood' in Burmese [as] something 
not uttered but implied; they just aren't there\dots
\medskip

AL Becker, {\bf Beyond Translation} \cite{3} ( p 8)}\end{quotation} \bigskip

{\bf 2.1} For the philologist Becker it is fundamental that languages (and 
cultures) have `exuberances and deficiencies' that can make translation 
problematic; in the present case these issues are acute. 

It seems generally agreed \cite{21} that the $\CF$ is underdetermined. One of
the first things a mathematician notices (\cf the original statement, 
included below as an appendix) is the absence of quantifiers, which usually 
convey information about the circumstances under which a proposition holds; 
moreover, the reversal of function and term [condition {\bf 2}] is mysterious. 
One issue - the assumption that the element $a$ has a natural or canonical 
dual $a^{-1}$ - can perhaps be resolved by reading the $\CF$ as 
{\it asserting}, in the context $\{x,a,y,b\}$, the {\it existence} of 
$a^{-1}$: for example, as in \cite{12}. On the other hand, an exuberance of 
the paraphrase above is a precise specification of the binary octahedral 
group as a repository for the values of our analog $\Phi$ of L\'evi-Strauss's 
$F$. 

{\bf 2.2} I make no claim about the uniqueness of the model proposed here:
I hope there are better ones.  Perhaps this is the place to say that
my concern with the $\CF$ is not its validity, or `truth-value'; it is rather
whether or not it can be usefully be interpreted as a formal mathematical
assertion. Its interest as an empirical hypothesis, like Bohr's model for 
the atom, seems well-established. \textcyr{Po moemu} the key question is what 
`interpretation', in this context, could even mean. 

{\bf 2.3} An answer may lie in the ancient opposition, older than Zeno, 
between continuous and discrete. An anthropologist, like William James' infant 
\cite{17} (Ch 13) is `assailed by eyes, ears, nose, skin and entrails at once, 
feels it all as one great blooming, buzzing confusion'; but `nowadays, 
fundamental psychological changes occur [to the mathematician]. Instead of 
sets, clouds of discrete elements, we envisage some sorts of vague spaces 
\dots mapped one to another. If you want a discrete set, then you pass to the 
set of connected components of spaces defined only up to homotopy' 
\cite{20} (p 1274), \cite{17}.

I believe these remarks of Manin point to an answer to this question, in
terms of cognitive condensation of chaotic clouds of experience into discrete, 
classifiable conceptual entities, \cf \: \cite{4}. \bigskip

{\bf \S 3} {\sc A short mathematical glossary}

This section summarizes more than enough background from the theory of groups
for the purposes of this paper:

{\bf 3.1.1} The objects of the category of groups (for example the
integers $\Z = \{\dots,-2,-1,0,1,2,\dots\}$) consist of sets $G$ with
an associated multiplication operation
\[
\mu_G : G \times G \ni (g,h) \mapsto g \cdot h \in G
\]
and an identity element $1 = 1_G \in G$, subject to familiar rules of 
associativity which I will omit. A homomorphism (or map, or morphism) 
$\phi : G \to G'$ between groups respects multiplication (\ie $\phi(g 
\cdot g') = \phi(g) \cdot \phi(g')$ etc.); a composition of homomorphisms 
is again a homomorphism, thus defining a category. The set of one-to-one
self-maps of a set with $n$ elements, for example, defines the symmetric
group $\Sigma_n$.

A group $A$ is commutative if $ a \cdot a' = a' \cdot a$ for any $a,a' \in A$;
in such cases the multiplication operation is often written additively, \ie
with $a + a'$ instead of $a \cdot a'$, \eg as in the additive group $\Z$ of
integers. The order of a group is the (not necessarily finite) number of its
elements.

{\bf Example} The set of isomorphisms $\alpha: G \to G$ with itself is 
similarly a group $\aut(G)$ (of automorphisms of $G$). Any element $h \in G$ 
defines an {\it inner} automorphism
\[
\alpha_h : G \ni g \mapsto hgh^{-1} \in G
\]
of $G$; it is a group homomorphism since
\[
\alpha_h(g \cdot g') = h(g \cdot g')h^{-1} = hgh^{-1} \cdot hg'h^{-1} =
\alpha_h(g) \cdot \alpha_h(g') \;,
\]
and $h \mapsto \alpha_h : G \to \aut(G)$ is a homomorphism since
\[
(\alpha_h \circ \alpha_k)(g) = \alpha_h(kgk^{-1}) = hkgk^{-1}h^{-1} =
hk \cdot g \cdot (hk)^{-1} = \alpha_{hk}(g) \;;
\]
much mathematics consists of shuffling parentheses. The subgroup $\iN(G) = 
\{ \alpha_g \in \aut(G) \}$ is {\it normal} in that 
\[
(\forall \beta \in \aut(G)) \; \alpha \in \iN(G) \Rightarrow \beta \alpha 
\beta^{-1} \in \iN(G) \;.
\]
In particular, if $G=A$ is commutative, then $\iN(A) = \{1\}$ is the trivial 
group.

If a subgroup $H$ of $G$ is normal, then the set $G/H = \{gH \;|\; g \in G \}$
(of `orbits' of elements of $G$ under right multiplication by $H$) is again
a group.\bigskip

{\bf 3.1.2} A composition
\[
\xymatrix{
H \ar[r]^\phi & G \ar[r]^\psi & K }
\]
of homomorphisms is {\it exact} if the {\it image}
\[
\im \; \phi = \{g \in G \:|\: \exists h \in H, \; \phi(h)= g \}
\]
of $\phi$ equals the {\it kernel}
\[
\ker \; \psi = \{g \in G \:|\: \psi(g) = 1_K \}
\]
of $\psi$; this implies in particular that the composition $\psi \circ \phi$
is trivial (\ie maps every element of $H$ to the identity element of $K$), but
is more restrictive. A sequence
\[
\xymatrix{
1 \ar[r] & H \ar[r]^\phi & G \ar[r]^\psi & K \ar[r] & 1 }
\]
of groups and homomorphisms is exact if its consecutive two-term compositions
are exact; this implies that 

$\bullet$ $\phi$ is one-to-one (or is a monomorphism, or has trivial kernel),

$\bullet$ $\psi$ is `onto' (or surjective, or $K = \im\; \psi$), and

$\bullet$ $\im \; H$ is normal in $G$, and $\psi$ factors through an 
isomorphism $G/H \cong K$. 

{\bf 3.1.3} Such an exact sequence is said to {\it split}, if there is a 
homomorphism $\rho : K \to G$ inverse to $\psi$ in the sense that the 
composition $\psi \circ \rho$
\[
\xymatrix{ 
K \ar[r]^\rho & G \ar[r]^\psi & K }
\]
is the identity map $K$. In that case there is a unique homomorphism
\[
\ve : K \to \aut(H)
\]
such that (the `semi-direct product) $G \cong H \rtimes K$ is isomorphic to 
the group defined on the set product $H \times K$, with twisted multiplication
\[
(h_0,k_0) \cdot (h_1,k_1) = (h_0 \cdot \ve(k_0)h_1,k_0k_1) \;;
\] 
such a split sequence will usually be displayed below as
\[
\xymatrix{
1 \ar[r] & H \ar[r] & G \cong H \rtimes K \ar[r] & K \ar[r] & 1 \;.}
\]

{\bf 3.2} {\sc Some topological groups and algebras}

{\bf 3.2.1} The quaternion group
\[
Q = \{\pm 1,\pm i,\pm j, \pm k \}
\]
of order eight is defined by three elements $i,j,k$ with multiplication $i^2 = 
j^2 =k^2 = -1$ and
\[
ij = k = -ji, \; jk = i = -kh, \; ki = j = -ik \;.
\]
The (noncommutative) division algebra
\[
\mH = \{q = q_0 + q_1i + q_2j + q_3k \:|\: q_i \in \R \} \cong \R^4
\]
of Hamiltonian quaternions is the four-dimensional real vector space 
with multiplication extended from $Q$; alternately, it is the two-dimensional 
complex vector space
\[
\mH = \{z_0 + z_1j \:|\: z_i \in \C\} \;,
\]
where $z_0 = q_0 + iq_1, \; z_1 = q_2 + iq_3$. The quaternions thus extend 
the field $\C$ of complex numbers much as $\C$ extends the field $\R$ of 
real numbers. 

The quaternion conjugate $q^* = q_0 - q_1i - q_2j - q_3k$ to $q$ has 
positive product
\[
q^* \cdot q = q \cdot q^* = |q|^2 = \sum q_i^2 > 0
\]
with $q$ if $q \neq 0$, implying the existence of a multiplicative inverse 
$q^{-1} = |q|^{-2}q^*$. This defines an isomorphism
\[
\mH^\times \ni q \mapsto (|q|^{-1}q,|q|) \in S^3 \times \R^\times_+
\]
making the three-dimensional sphere $S^3$ a group under multiplication.
This notation is nonstandard but convenient; note that $q^{**} = q$ and 
that quaternion conjugation $* : q \mapsto q^*$ is an {\it anti}-homomorphism, 
\ie 
\[
(u \cdot v)^* = v^* \cdot u^* \;.
\]
[Similarly,
\[
\R^\times \cong S^0 \times \R^\times_+, \; S^0 = \{\pm 1\},
\]
while
\[
\C^\times \cong S^1 \times \R^\times_+
\]
(where $S^n = \{ {\bf x} \in \R^{n+1} \:|\: |{\bf x}|^2 = 1 \}$ is the 
$n$-dimensional sphere of radius one, \eg the circle when $n=1$).] 

{\bf 3.2.2} The subalgebra (\ie closed under addition and multiplication, but 
not division) of {\it Lipschitz} quaternions in $\mH$ is the set of $q$ with 
integral coordinates ($q_i \in \Z$), while the subalgebra of {\it Hurwitz}
quaternions consists of elements $q$ with all coordinates either integral
or half-integral (\ie such that each $q_i$ is half of an {\it odd} integer). 
Finally, the subalgebra of Lipschitz (integral) quaternions has an additional
(commutative but nonassociative) Jordan algebra product
\[
u,v \mapsto \{u,v\} = \half (u \cdot v + v \cdot u) = \{v,u\}\;,
\]
\eg
\[
\{1,1\} = 1, \; \{i,i\} = \{j,j\} = \{k,k\} = -1, \; \{i,j\} = \{j,k\} = 
\{k,i\} = 0 \;.
\]
This allows us to define a (non-commutative) Jordan ratio 
\[
\{u:v\} = \{u,v^*\} = |v|^{-2}\{u,v^{-1}\}
\]
for $u,v \in \mH$, distributive
\[
\{u + u':v\} =  \{u:v\} + \{u:v'\}, \; \{u:v + v'\} = \{u:v\} + \{u:v'\}
\]
in both variables, satisfying 
\[
\{u:v\}^* = \{u^*:v^*\} \;.
\] 

{\bf Remark} $\mH$ can be regarded as a subalgebra of the $2 \times 2$ complex
matrices $M_2(\C)$, in such a way that the quaternion norm $|q|^2$ equals 
the determinant of $q$, regarded as a matrix. This identifies the 3-sphere
$S^3$ with a subgroup of the Lie group $\Sl_2(\C)$ of complex $2 \times 2$
matrices with determinant one; as such, it is the maximal compact (special 
unitary) subgroup $\SU(2)$ of $\mH^\times$. The special orthogonal group 
$\SO(3) \cong \SU(2)/\{\pm 1\}$ (of rotations in three dimensions) is a 
quotient of this group, and the various `binary' (tetrahedral, octahedral
etc.) groups lift the symmetry groups of the classical Platonic solids 
\cite{12} to subgroups of the three-sphere. See \cite{2}, and its comments, 
for some very pretty animations of a certain `24-cell' associated \cite{10} 
(\cf \: note {\it ar}) to the octahedral and tetrahedral groups. [The 
noncompact group $SL_2(\C)$ is similarly a double cover of (the identity 
component of) the physicists' Lorentz group.]

{\bf 3.3} The subset $Q \subset \mH^\times$ (of Lipschitz units) is a finite 
subgroup of $\SU(2)$. Similarly, the subset
\[
Q \subset A_{24} \subset \SU(2)
\]
(of {\it Hurwitz} units) is the union of $Q$ with the set of sixteen elements
of the form
\[
\half [\pm 1 \pm i \pm j \pm k] \;.
\]
It is known as well \cite{5} as the binary tetrahedral group $2 \cdot T$. 
\bigskip

Klein's (commutative) `Vierergruppe'
\[
V = \{1,I,J,K\}
\]
with multiplication $I^2 = J^2 = K^2 = 1$ and $IJ = JI = K, \; JK= KJ = I, \;    
KI = IK = J$ can be regarded as the subgroup
\[
V = \iN(Q) \subset \aut(Q)
\]
defined by $\alpha_i = I, \; \alpha_j = J, \; \alpha_k = K$.

{\bf Exercise:}

{\it i)} $I$ sends $i$ to itself, and $j,k$ to $-j,-k$; similarly $J$ sends
$j$ to itself, while $i,k \mapsto -i,-k$, etc. The cyclic permutation
\[
(abc) = (a \to b \to c \to a) \in \Sigma_3
\]
of three things defines a homomorphism $C_3 \cong \{1,\sigma,\sigma^2\} \to
\aut(Q)$ sending  $\sigma$ to $(ijk)$.

{\it ii)} The map $C_2^2 = C_2 \times C_2 \to V$ defined by $(1,0) \mapsto 
I, \; (0,1) \mapsto J, \; (1,1) \mapsto K$ (and $1 \mapsto (0,0)$) is a
(nonunique) isomorphism.

{\it iii)} For example, in $V \rtimes C_3 = A_4$ (\ie the alternating 
subgroup of order 12 (see below) of $\Sigma_4 \cong \aut(Q)$) we have  
\[
(I\sigma^2) \cdot (I\sigma^2) = I \sigma^2(I) \cdot \sigma^4 = IK\sigma = 
J\sigma \;.
\]
The anti-automorphism $\lambda = *I\sigma^2$ of $Q$ satisfies
\[
\lambda(i) = k, \; \lambda(j) = - i, \; \lambda(k) = j,
\]
\eg $\lambda(j) = (I\sigma^2(j))^* = (Ii)^* = i^* = - i$, \cf [21 \S 5]. 

{\bf 3.4.1} {\sc Some useful small groups} 

{\bf order {\; \;} \& {\; \:} name}

{\bf n} {\; \; \; \;} cyclic $C_n = \{0,\dots,n-1\}, \; n \geq 1$
\[
\xymatrix{
1 \ar[r] & \Z \ar[r]^n & \Z \ar[r] & \Z/n\Z \cong C_n \ar[r] & 1}
\]

{\bf n!} {\; \; \;} symmetric $\Sigma_n, \; n \geq 1$

{\bf 4} {\; \; \; \;}  Klein Vierergruppe $V$
\[
V = \{1,I,J,K\} \cong C_2 \times C_2 \; ({\rm non-uniquely})
\]

{\bf 6} {\; \; \; \;} symmetric $\Sigma_3$
\[
\xymatrix{
1 \ar[r] & C_3 \ar[r] & \Sigma_3 \cong C_3 \rtimes C_2 \ar[r] & C_2 \ar[r] & 1}
\]

{\bf 8} {\; \; \; \;}  (Lipschitz) quaternion units $Q = \{\pm 1,\pm i,\pm j,\pm k \}$
\[
\xymatrix{
1 \ar[r] & C_2 \ar[r] & Q \ar[r] & V \ar[r] & 1}
\]
% \[
% \xymatrix{
% 1 \ar[r] & C_4 \cong \{\pm 1, \pm i\} \ar[r] & Q \ar[r] & C_2 \ar[r] 1 }
% \]
{\bf 12} {\; \; \;}  alternating or tetrahedral ($T = A_4$)
\[
\xymatrix{
1 \ar[r] & V \ar[d]^= \ar[r] & A_4 \cong V \rtimes C_3 \ar[d] \ar[r] & C_3
\ar[d] \ar[r] & 1\\
1 \ar[r] & V \ar[r] & \Sigma_4 \cong V \rtimes \Sigma_3 \ar[r] & \Sigma_3
\ar[r] & 1}
\]

{\bf 24} {\; \; \;}  $\Sigma_4$ as above; binary tetrahedral $2 \cdot T$ =
Hurwitz units $A_{24}$
\[
\xymatrix{
1 \ar[r] & C_2 \ar[r] & 2 \cdot T = Q \rtimes C_3 \ar[r] & A_4 = V \rtimes C_3
\ar[r] & 1}
\]

{\bf 48} {\; \; \;} binary octahedral $2 \cdot O$
\[
\xymatrix{
1 \ar[r] & C_2 \ar[r] & 2 \cdot O  \ar[r] & 2 \cdot T \ar[r] & 1 }
\]
as well as
\[
\xymatrix{
1 \ar[r] & Q \ar[r] & 2 \cdot O \ar[r] & \Sigma_3 \ar[r] & 1 \;.}
\]

{\bf 3.4.2} The binary octahedral group $2 \cdot O$, regarded as a subgroup 
of the unit quaternions \cite{2,6}, is the disjoint union of $A_{24}$ with 
the set of twenty-four special elements 
\[
q = 2^{-1/2}[q_0 + q_1i + q_2j + q_3k]
\]
in which exactly two of $q_0,\dots,q_3$ are nonzero and equal $\pm 1$, 
in which $\Phi$ takes its values. 

{\bf 3.4.3} We have
\[
\aut(Q) \cong \Sigma_4 \cong \aut(2 \cdot T)
\]
and
\[
\aut(2 \cdot O) \cong C_2 \times (2 \cdot T \cong A_{24}) \;.
\]

{\bf 3.5} Some of the groups above can be presented in terms of matrices 
over finite Galois fields $\F$: in particular, $\Sigma_3 \cong \Sl_2(\F_2)$
and $A_{24} \cong \Sl_2(\F_3)$. Similarly, the binary icosahedral group
(which plays no role in this paper) is isomorphic to $\Sl_2(\F_5)$.

It is worth mentioning that the group of $2 \times 2$ matrices with entries
from $\Z$ and determinant one is a quotient
\[
\xymatrix{
1 \ar[r] & \Z \ar[r]^t & \B_3 \ar[r] & \Sl_2(\Z) \ar[r] & 1 }
\]
of Artin's three-strand braid group \cite{7}, which thus maps (by $\Z \to 
\Z/3\Z \cong \F_3$) to $\aut(2 \cdot O)$. 

[The set of braids on $n$ strands, imagined for example as displayed on a 
loom, define a group under `concatenation':

Technically, an $n$-strand braid can be defined as a smooth path in the space
of configurations defined by $n$ distinct points in the plane, starting for
example at time $t = 0$ at the integral points $(1,0),(2,0),\dots,(n,0)$
and ending at time $t = 1$ at the points $(1,1),(2,1),\dots,(n,1)$, though
not necessarily in that order. Such braids can be composed by concatenation
(\ie glueing and rescaling), and define elements of $\Sigma_n$ (sending $k$ 
to $l$ if the strand starting at $(k,0)$ ends at $(l,1)$). The braid group 
$\B_n$ is the set of such things under the equivalence relation roughly 
described as straightening: thus for example any braid can be parsed into a 
composition of elementary moves, in which one strand over- or under-passes 
one of its nearest neighbors.  For example, $\B_3$ can be presented as the 
group with two generators $a,b$ satisfying the `braid relation' $aba = bab$, 
thus the map $t$ above sends 1 to the full twist $(aba)^2 = (ab)^3$.]\bigskip

{\bf \S 4} {\sc A Calculation}\bigskip

{\bf 4.1} L\'evi-Strauss's formula is expressed in terms of formal analogies, 
\eg $F_x(a):F_y(b)$, understood roughly as a ratio, in a sense going back
to Eudoxus; but noncommutative algebra distinguishes the left fraction
$a^{-1}b = a\backslash b$ from the right fraction $ba^{-1} = b/a$. The
noncommutative ratio [3.3.2] splits this difference. Thus if $x,a,y,b \in Q$ 
we have

\newpage

\[
\{\Phi_x(a):\Phi_y(b)\} = \half \{x-a:y-b\} = \half [(\{x:y\} + \{a:b\}) - 
(\{a:y\} + \{x:b\})] \;,
\]
so for example if $x=1,\; a=i,\; y=j,\; b = k$ we have
\[
\{\Phi_x(a):\Phi_y(b)\} = \{\Phi_1(i):\Phi_j(k)\} = \half \{1-i:j-k\} = 
\]
\[
\half [(\{1:j\} + \{i:k\}) - (\{i:j\} + \{1:k\})] = \half [(j+0)-(0+k)]
= \half (j-k) \;,
\]
while 
\[
\{\Phi_x(b):\Phi_{a^{-1}}(y)\} = \half \{1-k:-i-j\} = 
\]
\[
\half [ (-\{1:i\} - \{1:j\}) + (\{k:i\} + \{k:j\})] = \half [-i - j + 0] 
= - \half (i+j) \;.
\]

But now, applying the anti-automorphisms $\lambda$ of \cite{20} (\S 5), we have
\[
\lambda(j - k) = - (i + j)
\]
by 3.3{\it iii} above, and the proposition is verified. $\Box$ \bigskip

{\bf Remarks}

{\bf 4.2.1} The automorphism group $\Sigma_4$ of $Q$ preserves the set of 
special elements [3.4.2] of $2 \cdot O$, \ie the possible values of $\Phi$, 
as well as their Jordan ratios; but it differs from the (inner) automorphism
group $A_{24}$ of $2 \cdot O$. However, the term $x$ appears in the $\CF$ in
the same place on both sides of the equation, and can be interpreted as 
playing the role of $1 \in Q$, fixed by all automorphisms. The canonical 
formula appears in variant forms in the literature, but (as far as I know) 
they all feature a term in the same place on both sides of the equation, 
allowing the variants to be reconciled by a cyclic permutation in $\Sigma_3 
\subset \Sigma_4$, which lies in a quotient of $\aut(2 \cdot O)$.  
  
{\bf 4.2.2} The classic work of Thom on singularity theory \cite{24} has 
turned out (\eg under the influence of Arnol'd, McKay, and others) to have 
deep connections with the theory of Platonic symmetry groups. The binary 
orthogonal group, in particular, seems to be related to a certain `symbolic 
umbilic' singularity \cite{9,11}, a special case of Thom's original 
classification.

{\bf 4.2.3} As a closing remark: the model proposed here is in fact not
that complicated. To a mathematician, perhaps the most interesting implication
is its connection with the theory of braids, which is arguably related to 
the processing of recursion, and to cognitive evolution.

{\bf Acknowledgements and thanis}: I am deeply indebted to many people, 
especially John Baez, Ellen Contini-Morava, Fred Damon, John McKay , Tony 
Phillips, Emily Riehl, Dale Rolfsen, Lucien Scubla, and Roy Wagner, for 
conversations, insight, and advice, during the preparation of this paper. 
Its excesses and deficiencies, though, are my own responsibilty. \newpage

{\bf Appendix} From {\sc The structural study of myth} \cite{14}, Journal of 
American Folklore 68 (1955):

7.30 Finally, when we have succeeded in organizing a whole series of
variants in a kind of permutation group, we are in a position to formulate
the law of that group. Although it is not possible at the present stage to
come closer than an approximate formulation which will certainly need to
be made more accurate in the future, it seems that every myth (considered
as the collection of all its variants) corresponds to a formula of the
following type:
\[
F_x(a):F_y(b) \; \simeq \; F_x(b):F_{a^{-1}}(y)
\]
where, two terms being given as well as two functions of these terms, it
is stated that a relation of equivalence still exists between two situations 
when terms and relations are inverted, under two conditions: {\bf 1.} that 
one term be replaced by its contrary; {\bf 2.} that an inversion be made 
between the {\it function} and the {\it term} value of the two elements. 
\bigskip

\bibliographystyle{amsplain}

\end{document}